\documentclass[12pt]{amsart}       
\usepackage{txfonts}
\usepackage{amssymb}
\usepackage{eucal}
\usepackage{bbm}
\usepackage{graphicx}
\usepackage{amsmath}
\usepackage{amscd}
\usepackage[all]{xy}           
\usepackage{amsfonts,latexsym}
\usepackage{xspace}
\usepackage{epsfig}
\usepackage{float}
\usepackage{mathrsfs}
\usepackage{amsbsy}
\usepackage{color}
\usepackage{shuffle}
\usepackage{fancybox}
\usepackage{colordvi}
\usepackage{multicol}
\usepackage{colordvi}
\usepackage{ifpdf}
\ifpdf
  \usepackage[colorlinks,final,backref=page,hyperindex]{hyperref}
\else
  \usepackage[colorlinks,final,backref=page,hyperindex,hypertex]{hyperref}
\fi
\usepackage[active]{srcltx} 

\usepackage{tikz}

\usepackage{graphicx}

\newtheorem{theorem}{Theorem}[section]

\newtheorem{prop-def}{Proposition-Definition}[section]
\newtheorem{coro-def}{Corollary-Definition}[section]

\theoremstyle{definition}
\newtheorem{definition}[theorem]{Definition}
\newtheorem{remark}[theorem]{Remark}

\newcommand{\nc}{\newcommand}
\nc{\tred}[1]{\textcolor{red}{#1}}
\nc{\tblue}[1]{\textcolor{blue}{#1}}
\nc{\tgreen}[1]{\textcolor{green}{#1}}
\nc{\tpurple}[1]{\textcolor{purple}{#1}}
\nc{\btred}[1]{\textcolor{red}{\bf #1}}
\nc{\btblue}[1]{\textcolor{blue}{\bf #1}}
\nc{\btgreen}[1]{\textcolor{green}{\bf #1}}
\nc{\btpurple}[1]{\textcolor{purple}{\bf #1}}
\nc{\NN}{{\mathbb N}}
\nc{\ncsha}{{\mbox{\cyr X}^{\mathrm NC}}} \nc{\ncshao}{{\mbox{\cyr
X}^{\mathrm NC}_0}}

\newcommand{\delete}[1]{}

\nc{\mlabel}[1]{\label{#1}}
\nc{\mcite}[1]{\cite{#1}}
\nc{\mref}[1]{\ref{#1}}
\nc{\meqref}[1]{\eqref{#1}}
\nc{\mbibitem}[1]{\bibitem{#1}}

\delete{
\nc{\mlabel}[1]{\label{#1}{\hfill \hspace{1cm}{\bf{{\ }\hfill(#1)}}}}
\nc{\mcite}[1]{\cite{#1}{{\bf{{\ }(#1)}}}}
\nc{\mref}[1]{\ref{#1}{{\bf{{\ }(#1)}}}}
\nc{\meqref}[1]{\eqref{#1}{{\bf{{\ }(#1)}}}}
\nc{\mbibitem}[1]{\bibitem[\bf #1]{#1}}
}
\nc{\sha}{{\mbox{\cyr X}}}  
\newfont{\scyr}{wncyr10 scaled 550}
\nc{\ssha}{\mbox{\bf \scyr X}}
\nc{\shap}{{\mbox{\cyrs X}}} 
\nc{\shpr}{\diamond}    
\nc{\shp}{\ast} \nc{\shplus}{\shpr^+}
\nc{\shprc}{\shpr_c}    
\nc{\dep}{\mrm{dep}} \nc{\lc}{\lfloor} \nc{\rc}{\rfloor}
\nc{\db}{\leq_{\rm db}} \nc{\bfk}{{\bf k}}

\nc{\cala}{{\mathcal A}} \nc{\calb}{{\mathcal B}}
\nc{\calc}{{\mathcal C}}
\nc{\cald}{{\mathcal D}} \nc{\cale}{{\mathcal E}}
\nc{\calf}{{\mathcal F}} \nc{\calg}{{\mathcal G}}
\nc{\calh}{{\mathcal H}} \nc{\cali}{{\mathcal I}}
\nc{\call}{{\mathcal L}} \nc{\calm}{{\mathcal M}}
\nc{\caln}{{\mathcal N}} \nc{\calo}{{\mathcal O}}
\nc{\calp}{{\mathcal P}} \nc{\calr}{{\mathcal R}}
\nc{\cals}{{\mathcal S}} \nc{\calt}{{\mathcal T}}
\nc{\calu}{{\mathcal U}} \nc{\calw}{{\mathcal W}} \nc{\calk}{{\mathcal K}}
\nc{\calx}{{\mathcal X}} \nc{\CA}{\mathcal{A}}

\nc{\fraka}{{\mathfrak a}} \nc{\frakA}{{\mathfrak A}}
\nc{\frakb}{{\mathfrak b}} \nc{\frakB}{{\mathfrak B}}
\nc{\frakc}{{\mathfrak c}}
\nc{\frakD}{{\mathfrak D}} \nc{\frakF}{\mathfrak{F}}
\nc{\frakf}{{\mathfrak f}} \nc{\frakg}{{\mathfrak g}}
\nc{\frakH}{{\mathfrak H}} \nc{\frakL}{{\mathfrak L}}
\nc{\frakM}{{\mathfrak M}} \nc{\bfrakM}{\overline{\frakM}}
\nc{\frakm}{{\mathfrak m}} \nc{\frakP}{{\mathfrak P}}
\nc{\frakN}{{\mathfrak N}} \nc{\frakp}{{\mathfrak p}}
\nc{\frakS}{{\mathfrak S}} \nc{\frakT}{\mathfrak{T}}
\nc{\frakX}{{\mathfrak X}}

\font\cyr=wncyr10 \font\cyrs=wncyr7
\nc{\li}[1]{\textcolor{red}{#1}}
\nc{\lir}[1]{\textcolor{red}{Li:#1}}
\nc{\yi}[1]{\textcolor{blue}{Yi: #1}}
\nc{\xing}[1]{\textcolor{purple}{Xing:#1}}
\nc{\revise}[1]{\textcolor{red}{#1}}
\nc{\nan}[1]{\textcolor{blue}{Nan:#1}}

\numberwithin{equation}{section}
\nc{\X}{{\bf X}}  \nc{\Y}{{\bf Y}} \nc{\Z}{{\bf Z}} \nc{\G}{{\bf G}} \nc{\U}{{\bf U}} \nc{\HA}{{\bf H}}
\nc{\x}{\mathbb{X}}
\nc{\s}{\mathbb{S}}
\nc{\al}{\alpha}
\nc{\rrp}{{\mathscr{C}}^{\alpha}_{{\rm red}}([0, T], V)}
\nc{\rrpp}{{\mathscr{C}}^{\beta}([0, T], V)}
\nc{\fan}[1]{\|#1\|}
\nc{\bfone}{{\bf 1}}

\nc{\crrp}{{\mathscr{D}}^{\alpha}_{\X, {\rm red}}([0, T], W)}
\nc{\crrrp}{{\tilde{\mathscr{D}}}^{\alpha}_{\X}([0, T], W)}

\nc{\crrpp}{{\mathscr{D}}^{\beta}_{\X}([0, T], W)}
\nc{\crrpx}{{\mathscr{D}}^{\alpha}_{\X}([0, T], W)}
\nc{\crrpa}{{\mathscr{D}}^{\alpha}_{\X}([0, T], \mathcal{L}(V, U))}
\nc{\crrrpa}{{\mathscr{D}}^{\alpha}_{\X}([0, T], \mathcal{L}(V, U))}

\nc{\W}{{\mathcal{Y}}}
\nc{\sym}{\operatorname{Sym}}
\nc{\anti}{\operatorname{Anti}}
\nc{\red}{\operatorname{red}}

\begin{document}

\title[Reduced rough paths and RDEs]{Rough differential equations and reduced rough paths: a Lie bracket characterization
}
%
\author{Nannan Li}
\address{School of Mathematics and Statistics, Lanzhou University
Lanzhou, 730000, China
}
\email{linn2024@lzu.edu.cn}

\author{Xing Gao$^{*}$}\thanks{*Corresponding author}
\address{School of Mathematics and Statistics, Lanzhou University
Lanzhou, 730000, China;
Gansu Provincial Research Center for Basic Disciplines of Mathematics
and Statistics, Lanzhou, 730070, China
}
\email{gaoxing@lzu.edu.cn}

\begin{abstract}
This paper studies rough differential equations from the viewpoint of reduced rough paths in the H\"older regime \(\frac13<\alpha\le\frac12\).  A reduced rough path retains the first level and the symmetric part of the second level, while discarding the antisymmetric L\'evy-area component.  We identify the precise obstruction to determining rough differential equation solutions from this reduced information. For an RDE driven by vector fields \(F_1,\ldots,F_d\), we prove that any two rough paths with the same reduced projection produce the same solution for every common initial value if and only if
\[
[F_i,F_j]=0, \quad 1\le i,j\le d.
\]
Thus the antisymmetric L\'evy area is irrelevant exactly in the commuting-vector-field case.
\end{abstract}

\makeatletter
\@namedef{subjclassname@2020}{\textup{2020} Mathematics Subject Classification}
\makeatother
\subjclass[2020]{
60L20, 
}
\keywords{rough path, reduced rough path, rough differential equation, Lie bracket}

\maketitle

\setcounter{section}{0}

\allowdisplaybreaks

\section{Introduction}

\subsection{Rough paths}
In 1998, T. Lyons introduced the theory of rough paths~\cite{Ly98} to provide a rigorous mathematical framework for solving differential equations driven by rough signals; see also~\cite{FH20, FV10, LQ02} for systematic developments. These are differential equations of the form
\begin{equation*}
dY_t = \sum_{i=1}^d F_i(Y_t) dX^i_t, \quad \forall t \in [0, T],
\end{equation*}
where $Y:[0,T]\to  \mathbb R^n$ is the solution path, $X:[0,T]\to \mathbb R^d$ is the driving path, and $F_i: \mathbb R^n \to  \mathbb R^n$ are smooth vector fields. Equivalently, one may write $F=(F_1,\ldots,F_d): \mathbb R^n\to\mathcal L(\mathbb R^d,  \mathbb R^n)$. In the H\"older regime \(0<\alpha\le\frac12\), the driving signal \(X\) is in general too irregular for classical integration, and one needs additional information beyond the path itself in order to define differential equations driven by \(X\).  Rough path theory, initiated by Lyons, provides such a framework by enhancing \(X\) to a rough path \(\mathbf X\), whose second-level component encodes the iterated integrals required for stable integration and for solving rough differential equations.  This theory has since been further developed in several directions~\cite{BG22, Gu04, Gu10, VRST25}.

Controlled rough paths, introduced by Gubinelli~\cite{Gu04}, provide the standard analytic framework for defining rough integrals and solving RDEs; see also~\cite{BG22,FH20,FV10}.

\subsection{Reduced rough paths and the RDE obstruction}
A rough path $(X, \mathbb{X})$ enriches the path $X$ by including both the symmetric and antisymmetric components of the second-order iterated integral:
\[
\mathbb{X}_{s,t} := \int_s^t X_{s,u} \otimes dX_u.
\]
This object satisfies the full Chen relation:
\[
\mathbb{X}_{s,t} - \mathbb{X}_{s,u} - \mathbb{X}_{u,t} = X_{s,u} \otimes X_{u,t}.
\]
In contrast, a reduced rough path $(X, \mathbb{S})$ (see Definition~\mref{def:rrp}) retains only the symmetric part of the second-level tensor~\cite{FH20}:
\[
\operatorname{Sym}(\mathbb{X}_{s,t}) := \frac{1}{2} \left( \mathbb{X}_{s,t} + \mathbb{X}_{s,t}^\top \right),
\]
which discards the antisymmetric (area-type) component:
\[
\operatorname{Anti}(\mathbb{X}_{s,t}) := \frac{1}{2} \left( \mathbb{X}_{s,t} - \mathbb{X}_{s,t}^\top \right).
\]
Here $\mathbb{X}_{s,t}^\top$ is the transpose of $\mathbb{X}_{s,t}$  in its matrix representation. The symmetric part $\operatorname{Sym}(\mathbb{X}_{s,t})$ captures the symmetric part of the iterated integral, i.e., information that does not depend on the order of integration indices, while the antisymmetric component $\operatorname{Anti}(\mathbb{X}_{s,t})$ encodes the antisymmetric or ``rotational'' component, also known as the L\'evy area.

The full enhancement of a path---which includes nontrivial second-order information such as L\'evy area---may contain more information than is needed in situations where only symmetric second-order information is relevant. Indeed, a symmetric second-order coefficient annihilates every antisymmetric tensor and therefore depends only on the symmetric part of the second level.

For rough differential equations, however, this symmetry is not automatic.  If the equation is
\[
 dY_t=\sum_{i=1}^d F_i(Y_t)dX_t^i,
\]
then the derivative of the integrand \(F(Y)\) is represented at the solution level by \(DF(Y)F(Y)\).  The antisymmetric part of this coefficient is governed by the Lie brackets of the vector fields.  Equivalently, in the RDE expansion, the antisymmetric part of the second level enters through the bracket terms.  Consequently, two full rough paths with the same reduced projection may lead to different RDE solutions whenever some bracket is nonzero.  The main purpose of this paper is to make this obstruction precise.

In this paper, we ask when the RDE solution is completely determined by the reduced rough path.  We prove a sharp characterization: this happens if and only if the driving vector fields commute.  The sufficiency of this condition is explained by the
symmetry of the second-order coefficient in the RDE expansion.  The essential point is the converse direction: if the solution is determined by reduced data for all driving rough paths and all initial values, then every Lie bracket must vanish.  This shows that the Lie-bracket contribution is exactly the obstruction to using reduced rough paths for general RDEs.\\

\noindent\textbf{Outline.}
This paper is structured as follows.
Section~\ref{sec:2} recalls reduced rough paths. Section~\ref{sec:3} contains the main result: a Lie-bracket characterization of when RDE solutions are determined by the reduced rough path (Theorem~\ref{thm:rde}).\\

\noindent\textbf{Notation.}
We work over the field $\mathbb{R}$ of real numbers, which serves as the base field for all vector spaces, tensor products, and linear maps. For a continuous path
\[
X: [0,T] \to V, \quad t \mapsto X_t,
\]
its increment over an interval $[s,t]$ is denoted by $X_{s,t} := X_t - X_s$.

\section{Reduced rough paths}\label{sec:2}
In this section, we recall the notion of a reduced rough path and its relation to the symmetric part of the second level of a full rough path.

 Let $V$ be a Banach space. For each $n\in \mathbb{Z}_{\ge 1}$, the space of symmetric tensors is denoted by
\[
\operatorname{Sym}(V^{\otimes n}) := \{ T \in V^{\otimes n} \mid \sigma \cdot T = T \text{ for all } \sigma \in S_n \},
\]
where $S_n$ is the symmetric group of degree $n$ and for $T = v_1 \otimes \cdots \otimes v_n\in V^{\otimes n}$,
$$\sigma \cdot T:= \sigma \cdot (v_1 \otimes \cdots \otimes v_n):=  v_{\sigma(1)} \otimes \cdots \otimes  v_{\sigma(n)} .$$
The {\bf symmetrization operator} is the linear map
\begin{equation*}
\operatorname{Sym} : V^{\otimes n} \to \operatorname{Sym}(V^{\otimes n}), \quad T\mapsto \frac{1}{n!} \sum_{\sigma \in S_n} \sigma \cdot T.
\end{equation*}

We use the following definition of a reduced rough path.

\begin{definition}~\cite{CP19, FH20}
Let $\alpha \in (1/3,1/2]$ and $\Delta_T:=\{(s,t)\mid 0\le s\le t\le T\}$. An $\al$-H\"{o}lder {\bf reduced rough path} is a pair $(S, \mathbb S)$, where
$$S:[0, T]\to V,\quad \mathbb S: \Delta_T\to \sym(V\otimes V),$$
satisfying:
\begin{enumerate}
\item Reduced Chen relation:
$$\mathbb S_{s,t}-\mathbb S_{s,u}-\mathbb S_{u,t}= \sym( S_{s,u}\otimes  S_{u,t}), \quad 0\le s\le u\le t\le T.$$
\item Analytic regularity:
$$\fan{S}_{\al}:=\sup_{0\le s<t\le T}{\|S_{s, t}\|\over |t-s|^{\al}}<\infty, \quad \fan{\mathbb S}_{2\al}:=\sup_{0\le s<t\le T}{\|\mathbb S_{s, t}\|\over |t-s|^{2\al}}<\infty.$$
\end{enumerate}
\mlabel{def:rrp}
\end{definition}

\section{Main result}\label{sec:3}
In this section, we identify the precise obstruction to determining RDE solutions from reduced data alone.

Let $W$ be a Banach space.  We take $V=\mathbb R^d$ with canonical basis $e_1,\ldots,e_d$, and write
\[
 F_i(y):=F(y)e_i,\quad 1\le i\le d.
\]
For two sufficiently regular vector fields $F_i,F_j$ on $W$, we use the convention
\begin{equation}
 [F_i,F_j](y):=DF_j(y)F_i(y)-DF_i(y)F_j(y).
\mlabel{eq:bracket}
\end{equation}

Let $\al \in(\frac{1}{3},\frac{1}{2}]$. An \(\alpha\)-H\"older rough path over \(\mathbb R^d\) is a pair
\[
\mathbf X=(X,\mathbb X),
\]
where \(X: [0,T]\to V\) and \(\mathbb X\colon\Delta_T\to V\otimes V\), such that
\[
\mathbb X_{s,t}=\mathbb X_{s,u}+\mathbb X_{u,t}+X_{s,u}\otimes X_{u,t},\quad 0\le s\le u\le t\le T,
\]
and $\|X\|_\alpha+\|\mathbb X\|_{2\alpha}<\infty$. Its reduced projection is defined by
\[
\operatorname{red}(\mathbf X):=\bigl(X,\operatorname{Sym}\mathbb X\bigr).
\]
By a solution of the RDE
\[
\mathrm dY_t=F(Y_t)\,\mathrm d\mathbf X_t, \quad Y_0=\xi,
\]
we mean a controlled rough path \(Y\) with prescribed Gubinelli derivative $Y'=F(Y)$ satisfying
\[
Y_t=\xi+\int_0^tF(Y_r)\,\mathrm d\mathbf X_r.
\]
For \(F\in C_b^3(W;\mathcal L(\mathbb R^d,W))\), this equation has a unique solution. We say that the  RDE associated with $F$ is \emph{determined by the reduced rough path} if, for every pair of rough paths $\mathbf X=(X,\mathbb X)$ and $\widetilde{\mathbf X}=(\widetilde X,\widetilde{\mathbb X})$ satisfying
\[
 \red(\mathbf X)=\red(\widetilde{\mathbf X}),
\]
and for every initial value $\xi\in W$, the corresponding RDE solutions coincide.

The following theorem is the main result.  It shows that reduced rough paths determine RDE solutions exactly in the commuting vector field case.  The nontrivial implication is the necessity: if the RDE solution is determined by the reduced rough path for all driving rough paths and all initial values, then the vector fields must commute.

\begin{theorem}\label{thm:rde}
Let $\alpha\in(\frac13,\frac12]$ and $F=(F_1,\ldots,F_d)\in C_b^3(W;\mathcal L(\mathbb R^d,W))$. The following statements
are equivalent.
\begin{enumerate}
\item\label{it:determined}
For every pair of \(\alpha\)-H\"older rough paths \(\mathbf X\) and \(\widetilde{\mathbf X}\) over \(\mathbb R^d\) satisfying
\[
\operatorname{red}(\mathbf X)
=
\operatorname{red}(\widetilde{\mathbf X}),
\]
and for every initial value \(\xi\in W\), the corresponding RDE solutions coincide.
\item\label{it:commuting}
The vector fields commute
\[
[F_i,F_j]=0, \quad 1\le i,j\le d.
\]
\end{enumerate}
\end{theorem}

\begin{proof}
For \(z\in W\), define
\[
B_z(v,w):=D(F(\cdot)w)(z)[F(z)v],\quad v,w\in\mathbb R^d.
\]
Then
\[
B_z(e_i,e_j)-B_z(e_j,e_i)=[F_i,F_j](z).
\]

Assume first that
\[
[F_i,F_j]=0,
\qquad 1\le i,j\le d.
\]
Then \(B_z\) is symmetric for every \(z\in W\).
Let
\[
\mathbf X=(X,\mathbb X),\quad \widetilde{\mathbf X}=(X,\widetilde{\mathbb X})
\]
have the same reduced projection. Let \(Y\) and \(\widetilde Y\) denote the RDE solutions driven by \(\mathbf X\) and \(\widetilde{\mathbf X}\), respectively, with the common initial value \(\xi\). The solution \(Y\) is controlled by \(X\) with Gubinelli derivative
\[
Y'=F(Y).
\]
Consequently, \(F(Y)\) is controlled by \(X\) with Gubinelli derivative \(B_Y\).  Recall that the standard rough integral of a controlled path \(G\), with Gubinelli derivative \(G'\), against \(\mathbf X=(X,\mathbb X)\) is given by
\[
\int_s^t G_r\,\mathrm d\mathbf X_r=\lim_{|\pi|\to0}\sum_{[u,v]\in\pi}\left(G_uX_{u,v}+G'_u\mathbb X_{u,v}\right),
\]
where the limit is taken over partitions \(\pi\) of \([s,t]\)~\cite{FH20}. Set
\[
A_{s,t}:=\widetilde{\mathbb X}_{s,t}-\mathbb X_{s,t}.
\]
Since the reduced projections agree, \(A_{s,t}\) is antisymmetric. Therefore, for every partition \(\pi\) of \([s,t]\), the symmetry of \(B_{Y_u}\) gives 
\begin{align*}
\sum_{[u,v]\in\pi}\left(F(Y_u)X_{u,v}+B_{Y_u}\widetilde{\mathbb X}_{u,v}\right)-\sum_{[u,v]\in\pi}\left(F(Y_u)X_{u,v}+ B_{Y_u}\mathbb X_{u,v}\right)=\sum_{[u,v]\in\pi}B_{Y_u}A_{u,v}=0.
\end{align*}
Passing to the limit yields
\[
\int_s^tF(Y_r)\,\mathrm d\widetilde{\mathbf X}_r=\int_s^tF(Y_r)\,\mathrm d\mathbf X_r.
\]
Thus \(Y\) also solves the RDE driven by \(\widetilde{\mathbf X}\). Uniqueness for the RDE driven by \(\widetilde{\mathbf X}\) implies that $Y=\widetilde Y$ on \([0,T]\).  This proves (\ref{it:commuting}) \(\Rightarrow\) (\ref{it:determined}).

Conversely, assume (\ref{it:determined}) is true. Fix \(i\ne j\) and \(a\ne0\), and define
\[
A^{ij,a}_t\equiv0, \quad \mathbb A^{ij,a}_{s,t}:=a(t-s)(e_i\otimes e_j-e_j\otimes e_i).
\]
The second level is additive, so
\[
\mathbb A^{ij,a}_{s,t}=\mathbb A^{ij,a}_{s,u}+\mathbb A^{ij,a}_{u,t},
\]
and hence Chen's relation holds because the first level vanishes. Moreover, since \(2\alpha\le1\),
\[
|\mathbb A^{ij,a}_{s,t}|\le C|a|T^{1-2\alpha}|t-s|^{2\alpha}.
\]
Thus
\[
\mathbf A^{ij,a}=(A^{ij,a},\mathbb A^{ij,a})
\]
is an \(\alpha\)-H\"older rough path, and
\[
\operatorname{red}(\mathbf A^{ij,a})=(0,0)=\operatorname{red}(\mathbf0).
\]
Let the common initial value be \(y\in W\). The solution driven by \(\mathbf0\) is the constant path
\[
Y_t\equiv y.
\]
By (\ref{it:determined}), the solution driven by \(\mathbf A^{ij,a}\) is the same constant path. In the controlled path formulation its Gubinelli derivative is
\[
Y'=F(y),
\]
and the derivative of the integrand \(F(Y)\) is \(B_y\). Therefore, for every \(t>0\),
\begin{align*}
0=Y_t-y&=\int_0^tF(y)\,\mathrm d\mathbf A^{ij,a}_r\\
&=\lim_{|\pi|\to0}\sum_{[u,v]\in\pi}\left(F(y)A^{ij,a}_{u,v}+B_y\mathbb A^{ij,a}_{u,v}\right)\\
&=\lim_{|\pi|\to0}\sum_{[u,v]\in\pi}B_y\,\mathbb A^{ij,a}_{u,v}\\
&=at\bigl(B_y(e_i,e_j)-B_y(e_j,e_i)\bigr)\\
&=at[F_i,F_j](y).
\end{align*}
Since \(a\ne0\), \(t>0\), and \(y\in W\) is arbitrary,
\[
[F_i,F_j]=0.
\]
This proves
(\ref{it:determined}) \(\Rightarrow\) (\ref{it:commuting}).
\end{proof}

\begin{remark}
The proof also explains the nontrivial direction of the theorem.  The sufficiency follows from the fact that, when the brackets vanish, the second-order coefficient \(B_z\) is symmetric and therefore cannot see the antisymmetric part of \(\mathbb X\).  Conversely, the pure-area rough paths \(\mathbf A^{ij,a}\) have the same reduced projection as the trivial rough path but produce the exact contribution \(at[F_i,F_j](y)\).  Hence any nonzero Lie bracket is detected by the antisymmetric second level and prevents the RDE solution from being determined by the reduced rough path.
\end{remark}

\vskip 0.2in

\noindent
{\bf Acknowledgments.} This work is supported by the National Natural Science Foundation of China (No.~12571019), the Natural Science Foundation of Gansu Province (No.~25JRRA644), Innovative Fundamental Research Group Project of Gansu Province (No.~23JRRA684) and Longyuan Young Talents of Gansu Province.

\noindent
{\bf Declaration of interests. } The authors have no conflicts of interest to disclose.

\noindent
{\bf Data availability. } Data sharing is not applicable as no new data were created or analyzed.

\end{document}